\documentclass{llncs}

\usepackage[utf8]{inputenc}
\usepackage[english]{babel}

\usepackage{amsmath} 
\usepackage{amsbsy} 
\usepackage{amssymb}
\usepackage{bm}
\usepackage{cite}

\newcommand{\grad}{\mathop{\rm grad}\nolimits}
\renewcommand{\div}{\mathop{\rm div}\nolimits}
\newcommand{\const}{\mathop{\rm const}\nolimits}

\title{\bf Flux-splitting schemes for parabolic problems}

\author{Petr N. Vabishchevich}

\institute{Nuclear Safety Institute, 52, B. Tulskaya, 115191 Moscow, Russia}

\begin{document}

\maketitle

\begin{abstract}

To solve numerically boundary value problems for parabolic equations with mixed
derivatives, the construction of difference schemes with prescribed quality 
faces essential difficulties.
In parabolic problems, some possibilities are associated with the transition to
a new formulation of the problem, where the fluxes 
(derivatives with respect to a spatial direction)
are treated as unknown quantities. In this case,
the original problem is rewritten in the form of a boundary value problem for the system
of equations in the fluxes. This work deals with studying schemes 
with weights for parabolic
equations written in the flux coordinates. Unconditionally stable flux
locally one-dimensional schemes of the first and second order of approximation in time 
are constructed for parabolic equations without mixed
derivatives. A peculiarity of the system of equations written in flux variables
for equations with mixed derivatives is that there do exist coupled terms with
time derivatives.

\end{abstract}

\section{Introduction}

Investigating many applied problems, we can consider a second-order parabolic equation
with mixed derivatives as the basic equation. An example is diffusion processes
in anisotropic media. In desining various approximations
for the corresponding boundary-value problems, we focus on the inheritance
of the primary properties of the differential problem 
during the construction of the discrete problem.

Locally one-dimensional difference schemes are obtained in a simple enough way
for second-order parabolic equations without mixed derivatives
\cite{Marchuk,SamarskiiVabischevich1999}.
Mixed derivatives complicate essentially the construction of unconditionally stable sche\-mes
of splitting with respect to the spatial variables for parabolic equations with
variable derivatives, even for two-dimensional problems.

In some problems, it is convenient to use the fluxes 
(derivatives with respect to a spatial direction) as
unknow quantities. 
This idea may be implemented in the most simple manner for one-dimensional problems
\cite{DegFav1}. To introduce fluxes, mixed and hybrid finite elements are applied
\cite{BrezziFortin1991,RobertsThomas1991}.
The original parabolic equation with mixed derivatives may be written
as a system of equations for the fluxes.
The basic peculiarity of this system is that the time derivatives for the fluxes
in separate equations are interconnected to each other.
For the problem in the flux variables, unconditionally stable schemes with weights
are developed. Locally one-dimensional schemes are proposed
for problems without mixed derivatives.

\section{Differential problem}

In a bounded domain $\Omega$, the unknown function $u(\bm{x},t)$,
$\bm{x} = (x_1, x_2, ..., x_m)$,
satisfies the equation
\begin{equation}\label{2.1}
   \frac{\partial u}{\partial t} 
   - \sum_{\alpha,\beta  =1}^{m}
   \frac{\partial }{\partial x_\alpha} 
   \left ( k_{\alpha\beta  }({\bm x})  \frac{\partial u}{\partial x_\beta } \right ) = f({\bm x},t),
   \quad {\bm x}\in \Omega,
   \quad 0 < t \leq  T .
\end{equation}
Assume that the coefficients $k_{\alpha\beta}, \ \alpha,\beta = 1,2, ..., m$ satisfy
the conditions
\begin{equation}\label{2.2}
\begin{split}
  \underline{k} \sum_{\alpha =1}^{m} \xi^2_\alpha  ({\bm x}) \leq
  \sum_{\alpha,\beta  =1}^{m} k_{\alpha\beta  }({\bm x}) \xi_\alpha  ({\bm x}) \xi_\beta ({\bm x}) \leq
  \overline{k} \sum_{\alpha =1}^{m} \xi^2_\alpha  ({\bm x}), \\
  k_{\alpha\beta}=k_{\beta\alpha}, 
  \quad \alpha,\beta = 1,2, ..., m,
  \quad {\bm x} \in \Omega  
\end{split}
\end{equation}
for any $\xi_\alpha ({\bm x}), \ \alpha = 1,2, ..., m$ with constant $\underline{k} > 0$.
Consider the boundary value problem for equation (\ref{2.1}) 
with homogeneous Dirichlet boundary conditions
\begin{equation}\label{2.3}
   u({\bm x},t) = 0,
   \quad {\bm x}\in \partial \Omega,
   \quad 0 < t \leq  T
\end{equation}
and the initial conditions in the form
\begin{equation}\label{2.4}
   u({\bm x},0) = u^0({\bm x}),
   \quad {\bm x}\in \Omega.
\end{equation}

We introduce a vector quantity
$ \bm{q} = (q_1, q_2, ..., q_m)^T$ 
(the index $T$ denotes transposition) such that
\begin{equation}\label{2.5}
  \bm{q} = - \mathcal{K} \grad u,
\end{equation}
where $\mathcal{K} = (k_{\alpha\beta  })$ is a square matrix $m \times m$
($\mathcal{K} \in \mathbb{R}^{mm}$)
with elements $k_{\alpha\beta  }({\bm x}), \ \alpha,\beta = 1,2, ..., m$.
Using this notation, equation (\ref{2.1}) may be written as
\begin{equation}\label{2.6}
  \frac{\partial u}{\partial t} + \div \bm{q} = f,
  \quad \bm{x} \in \Omega,
  \quad 0 < t \leq T .
\end{equation}

We can write the above problem (\ref{2.3})--(\ref{2.5}) in the operator form.
Scalar functions are considered in the Hilbert space
$\mathcal{H} = L_2(\Omega)$ with the scalar product
and norm defined by the rules
\[
  (u,v) = \int_{\Omega } u(\bm{x}) v(\bm{x}) d \bm{x},
  \quad \|u\| = (u,u)^{1/2} . 
\] 
For vector functions, we use the Hilbert space
$\mathcal{V} = \mathbf{L}_2(\Omega)$, where
\[
  (\mathbf{q} ,\mathbf{g} ) = \sum_{\alpha =1}^{m} \int_{\Omega } q_\alpha (\bm{x}) g_\alpha (\bm{x}) d \bm{x},
  \quad \|\mathbf{q}\| = (\mathbf{q},\mathbf{q})^{1/2} . 
\] 

Taking into account (\ref{2.2}), we can treat the matrix $\mathcal{K}$ as a linear,
bounded, self-adjoint, and positive definite
operator in $\mathcal{V}$:
\begin{equation}\label{2.7}
  \mathcal{K}: \mathcal{V} \rightarrow \mathcal{V},
  \quad \mathcal{K} = \mathcal{K}^*,
  \quad \underline{k} \mathcal{E} \leq \mathcal{K} \leq \overline{k} \mathcal{E},
  \quad \underline{k} > 0, 
\end{equation}
where $\mathcal{E}$ is the identity operator in $\mathcal{V}$.
Suppose $\mathcal{D} u = - \grad u$,i.e.,
\begin{equation}\label{2.8}
  \mathcal{D}: \mathcal{H} \rightarrow \mathcal{V},
  \quad \mathcal{D} = \left ( - \frac{\partial }{\partial x_1} ,   
  - \frac{\partial }{\partial x_2} , \dots, - \frac{\partial }{\partial x_m} \right )^T .  
\end{equation}
On the set of functions that satisfy the boundary conditions (\ref{2.3}),
for the gradient and divergence operators, we have
\[
  \int_{\Omega}  u \div \mathbf{q}\ d \bm{x} +
  \int_{\Omega} \mathbf{q} \grad u \ d \bm{x} = 0 .
\]
It follows from this that $\mathcal{D}^* \mathbf{q}  =  \div  \mathbf{q} $, i.e.,
\begin{equation}\label{2.9}
  \mathcal{D}^*: \mathcal{V} \rightarrow \mathcal{H},
  \quad \mathcal{D}^* = \left (\frac{\partial }{\partial x_1} ,   
  \frac{\partial }{\partial x_2} , \dots, \frac{\partial }{\partial x_m} \right ) .  
\end{equation} 

In the above notation (\ref{2.7})--(\ref{2.9}),
from (\ref{2.3})--(\ref{2.5}), we obtain the Cauchy problem for the system
of operator-differential equations
\begin{equation}\label{2.10}
  \frac{d u }{d t} + \mathcal{D}^* \mathbf{q}  = f(t),
  \quad 0 < t \leq  T ,
\end{equation} 
\begin{equation}\label{2.11}
  \mathbf{q}  = \mathcal{K} \mathcal{D}  u,
\end{equation} 
\begin{equation}\label{2.12}
  u(0) = u^0 . 
\end{equation}
For the problem (\ref{2.1})--(\ref{2.4}), the following equation corresponds 
\begin{equation}\label{2.13}
  \frac{d u }{d t} + \mathcal{D}^* \mathcal{K} \mathcal{D}  u  = f(t),
  \quad 0 < t \leq  T ,
\end{equation}
wich is supplemented by the initial condition (\ref{2.12}).
Taking into account that
\[
 \frac{d}{d t} \mathcal{K} - \mathcal{K} \frac{d}{d t} = 0 , 
\]
it is possible to eliminate $u$ from the system of equations (\ref{2.10}), (\ref{2.11}) 
that gives
\begin{equation}\label{2.14}
  \mathcal{C} \frac{d \mathbf{q} }{d t} + \mathcal{D} \mathcal{D}^* \mathbf{q} = 
  \mathcal{D} f,
  \quad \mathcal{C} = \mathcal{K}^{-1},
  \quad 0 < t \leq  T .
\end{equation}
In view of (\ref{2.11}) and (\ref{2.12}), we put
\begin{equation}\label{2.15}
  \mathbf{q}(0) = \mathbf{q}^0 \equiv \mathcal{K} \mathcal{D}  u^0 .
\end{equation} 

In constructing locally one-dimensional schemes 
(schemes based on splitting with respect to spatial directions),
we focus on the coordinatewise formulation of equations
(\ref{2.10}), (\ref{2.11}), (\ref{2.14}) ) and (\ref{2.14}).
Let
\[
 \mathcal{D} = (\mathcal{D}_1, \mathcal{D}_2, \dots, \mathcal{D}_m)^T,
 \quad \mathcal{K}  = (\mathcal{K}_{\alpha\beta} ),
 \quad \mathcal{C}  = (\mathcal{C}_{\alpha\beta} ), 
\]
then the basic system of equations (\ref{2.10}), (\ref{2.11}) takes the form
\begin{equation}\label{2.16}
  \frac{d u }{d t} + \sum_{\alpha=1}^{m} \mathcal{D}_\alpha ^* q_\alpha   = f(t),
  \quad 0 < t \leq  T ,
\end{equation} 
\begin{equation}\label{2.17}
  q_\alpha   = \sum_{\beta=1}^{m} \mathcal{K}_{\alpha\beta}  \mathcal{D}_\beta   u,
  \quad \alpha = 1,2, ..., m .  
\end{equation}
The equation (\ref{2.13}) for $u$ is reduced to
\begin{equation}\label{2.18}
  \frac{d u }{d t} + \sum_{\alpha,\beta =1}^{m} \mathcal{D}_\alpha ^* \mathcal{K}_{\alpha\beta}  \mathcal{D}_\beta   u  = f(t),
  \quad 0 < t \leq  T .
\end{equation}
For the flux components (see (\ref{2.14})), we obtain
\begin{equation}\label{2.19}
  \sum_{\beta=1}^{m} \mathcal{C}_{\alpha\beta} \frac{d q_\beta  }{d t} + 
  \sum_{\beta=1}^{m}\mathcal{D}_\alpha  \mathcal{D}^*_\beta  q_\beta  = 
  \mathcal{D}_\alpha  f,
  \quad 0 < t \leq  T .
\end{equation} 

The equations of the system (\ref{2.19}) are connected with each other,
and, moreover, the time derivatives are interconnected.
The problem (\ref{2.12}), (\ref{2.18}) seems to be much easier --- 
we have a single equation instead of the system of $m$ equations.
Nevertheless, some possibilities to design locally one-dimensional schemes for
the system of equations are still there.

Here we present elementary a priori estimates for the solution
of the above Cauchy problems for operator-differential equations,
which will serve us as a checkpoint in the study of discrete problems.
Multiplying equation (\ref{2.13}) scalarly in $\mathcal{H}$ by $u$, we obtain
\[
  \|u \| \frac{d }{d t}\|u\| + 
  (\mathcal{K} \mathcal{D} u, \mathcal{D} u )
  = (f,u) .
\]
Taking into account (\ref{2.7}) and
\[
  (f,u) \leq   \|f \| \|u \|,
\]
we arrive at
\[
  \frac{d }{d t}\|u\| \leq  \|f \| .
\]
This inequality implies the estimate
\begin{equation}\label{2.20}
  \|u(t)\| \leq \|u^0\| + \int_{0}^{t}  \|f(\theta)\| d \theta  
\end{equation}
for the solution of the problem \ref{2.12}), (\ref{2.18}).

Now we investigate the problem (\ref{2.14}), (\ref{2.15}).
By the properties (\ref{2.7}) of the operator $\mathcal{K}$, 
for $\mathcal{C}$, we have
\begin{equation}\label{2.21}
  \mathcal{C}: \mathcal{V} \rightarrow \mathcal{V},
  \quad \mathcal{C} = \mathcal{C}^*,
  \quad \underline{c} \mathcal{E} \leq \mathcal{C} \leq \overline{c} \mathcal{E},
  \quad \underline{c} = \overline{k}^{\,-1} > 0,
  \quad \overline{c} = \underline{k}^{-1} . 
\end{equation}
In view of (\ref{2.21}), we define the Hilbert space $\mathcal{V}_\mathcal{C}$, 
where the scalar product and norm are
\[
  (\mathbf{q} ,\mathbf{g} )_\mathcal{C} = (\mathcal{C} \mathbf{q} ,\mathbf{g} ),
  \quad \|\mathbf{q}\|_\mathcal{C} = (\mathbf{q},\mathbf{q})_\mathcal{C}^{1/2} . 
\]
Multiplying equation (\ref{2.15}) scalarly in $\mathcal{V}$ by $\mathbf{q}$, we obtain 
\[
   \|\mathbf{q} \|_\mathcal{C} \frac{d }{d t}\|\mathbf{q}\|_\mathcal{C} + 
  (\mathcal{D}^* \mathbf{q}, \mathcal{D}^* \mathbf{q} )
  = (\mathcal{D}f,\mathbf{q}) .
\] 
In view of 
\[
 (\mathcal{D}f,\mathbf{q}) \leq  \|\mathcal{D}f\|_\mathcal{K} \|\mathbf{q}\|_\mathcal{C}, 
\]
we arrive at a priori estimate
\begin{equation}\label{2.22}
  \|\mathbf{q}(t)\|_\mathcal{C} \leq \|\mathbf{q}^0\| _\mathcal{C}+ \int_{0}^{t}  
  \|\mathcal{D}f(\theta)\|_\mathcal{K} \, d \theta  
\end{equation}
for the solution of the problem (\ref{2.14}), (\ref{2.15}).

\section{Approximation in space}

We conduct a detailed analysis using
a model two-dimensional parabolic problem in a rectangle
\[
  \Omega = \{ \bm{x} \ | \ \bm{x} = (x_1, x_2),
  \quad 0 < x_{\alpha} < l_{\alpha}, \quad \alpha = 1,2 \} .
\]
In $\Omega$, we introduce a uniform rectangular grid
\[
   \overline{\omega} = \{ \bm{x} \ | \ \bm{x} = (x_1, x_2),
   \quad x_\alpha = i_\alpha h_\alpha,
   \quad i_\alpha = 0,1,...,N_\alpha,
   \quad N_\alpha h_\alpha = l_\alpha\} 
\]
and let $\omega$ be the set of interior nodes
($\overline{\omega} = \omega \cup \partial \omega$). 
On this grid, scalar grid functions are given.
For grid functions $y(\bm{x}) = 0, \ \bm{x} \in \partial \omega$,
we define the Hilbert space $H = L_2({\omega})$
with the scalar product and norm
\[
  (y,w) \equiv \sum_{{\bm x} \in \omega}
  y({\bm x}) w({\bm x}) h_1 h_2,
  \quad \|y\| \equiv (y,y)^{1/2} .
\]

To determine vector grid functions, we have two main possibilities.
The first approach deals with specifying  vector functions on the same grid as
it used for scalar functions. The second possibility, which is traditionally widely
used, e.g., in computational fluid dynamics, is based on
the grid arrangement, where each individual component of a vector quantity is referred to its own mesh.
Here we restrict ourselves to the use of the same grid for all quantities, in particular,
for setting the coefficients $k_{\alpha\beta  }({\bm x}), \ \alpha,\beta = 1,2, ..., m$. 

Consider approximations for the differential operators
\[
   \mathcal{L}_{\alpha\beta  } u =
   - \frac{\partial }{\partial x_\alpha} 
   \left ( k_{\alpha\beta  }({\bm x})  \frac{\partial u}{\partial x_\beta } \right ) ,
   \quad \alpha,\beta = 1,2, ..., m .
\]
We apply the standard index-free notation from the theory of difference schemes
\cite{Samarskii1989} for the difference operators:
\[
  u_x = \frac{u(x+h) - u(x)}{h},
  \quad u_{\overline{x} } = \frac{u(x) - u(x-h)}{h} .
\]
If we set the coefficients of the elliptic operator at the grid points,
then
\begin{equation}\label{3.1}
  L_{\alpha\alpha} y = - \frac{1}{2} (k_{\alpha\alpha} u_{x_\alpha})_{\overline{x}_\alpha } -
  \frac{1}{2} (k_{\alpha\alpha} u_{\overline{x}_\alpha })_{x_\alpha} ,
  \quad \alpha =1,2 .  
\end{equation} \
More opportunities are available in approximation of operators with mixed
derivatives. As the basic discretization \cite{Samarskii1989}, we emphasize
\begin{equation}\label{3.2}
  L^{(1)}_{\alpha\beta  } y =  
  - \frac{1}{2} (k_{\alpha\beta} u_{x_\alpha})_{\overline{x}_\beta } -
  \frac{1}{2} (k_{\alpha\beta} u_{\overline{x}_\alpha })_{x_\beta} ,
\end{equation} 
\begin{equation}\label{3.3}
  L^{(2)}_{\alpha\beta  } y =  
  - \frac{1}{2} (k_{\alpha\beta} u_{x_\alpha})_{x_\beta } -
  \frac{1}{2} (k_{\alpha\beta} u_{\overline{x}_\alpha })_{\overline{x}_\beta} ,
  \quad \alpha, \beta  =1,2 , 
  \quad \alpha \neq  \beta .
\end{equation}
Instead of $L^{(1)}_{\alpha\beta}, L^{(2)}_{\alpha\beta}$, we can take their linear combination.
In particular, it is possible  \cite{matus_difference_2004} to put
\begin{equation}\label{3.4}
  L^{(3)}_{\alpha\beta  } = \frac{1}{2} L^{(1)}_{\alpha\beta  } + \frac{1}{2}  L^{(2)}_{\alpha\beta  } ,
  \quad \alpha, \beta  =1,2 , 
  \quad \alpha \neq  \beta . 
\end{equation}
In the general case, we set
\begin{equation}\label{3.5}
  L_{\alpha\beta  } = \chi  L^{(1)}_{\alpha\beta  } + (1-\chi)  L^{(2)}_{\alpha\beta  } ,
  \quad \alpha, \beta  =1,2 , 
  \quad \alpha \neq  \beta ,
  \quad \chi = \const.   
\end{equation}
The introduced discrete operators approximate the corresponding differential operators with the second order:
\begin{equation}\label{3.6}
  L_{\alpha\alpha} u = 
  \mathcal{L}_{\alpha\alpha} u + \mathcal{O} (h_{\alpha}^2),
  \quad L_{\alpha\beta  } u = \mathcal{L}_{\alpha\beta  } + \mathcal{O} (h^2),
  \quad \beta \neq \alpha, 
  \quad \alpha, \beta  =1,2 , 
\end{equation}
where $h^2 = h_1^2 + h_2^2$.

We define a grid subset $\overline{\omega}$, where
the corresponding components of vector quantities are defined. Let
\[
   \omega_1^+ = \{ \bm{x} \ |  
   \ x_1 = i_1 ,
   \ i_1 = 0,1,...,N_1-1,
   \ x_2 = i_2 h_2,
   \ i_2 = 1,2,...,N_2-1 \} ,
\]
\[
   \omega_1^- = \{ \bm{x} \ |  
   \ x_1 = i_1 ,
   \ i_1 = 1,2,...,N_1,
   \ x_2 = i_2 h_2,
   \ i_2 = 1,2,...,N_2-1 \} ,
\]
\[
   \omega_2^+ = \{ \bm{x} \ | 
   \ x_1 = i_1 h_1,
   \ i_1 = 1,2,...,N_1-1,
   \ x_2 = i_2 h_2,
   \ i_2 = 0,1,...,N_2-1 \}  ,
\]
\[
   \omega_2^- = \{ \bm{x} \ | 
   \ x_1 = i_1 h_1,
   \ i_1 = 1,2,...,N_1-1,
   \ x_2 = i_2 h_2,
   \ i_2 = 1,2,...,N_2 \}  ,
\]
and
\[
   \widetilde{\omega} = \omega_1^+ \cup  \omega_1^- \cup \omega_2^+ \cup \omega_2^- .
\]
For the grid vector variables, instead of two components, we will use four components, putting
\[
  \mathbf{q} = (q_1^+, q_1^-, q_2^+, q_2^-)^T,
  \quad q_\alpha^\pm =   q_\alpha^\pm({\bm x}), \quad {\bm x} \in \omega_\alpha^\pm ,
  \quad \alpha=1,2 .  
\] 

For the grid functions defined on grids
$\omega_{\alpha}^\pm, \ \alpha = 1,2$, we define the Hilbert spaces
$H_{\alpha}^\pm, \ \alpha = 1,2$, where
\[
  (y,w)_{\alpha}^\pm \equiv \sum_{{\bm x} \in \omega_{\alpha}^\pm}
  y({\bm x}) w({\bm x}) h_1 h_2,
  \quad \|y\|_{\alpha}^\pm \equiv ((y,y)_{\alpha}^\pm)^{1/2}, 
  \quad  \alpha = 1,2.
\]
For the grid vector functions in $V = H_1^+ \oplus H_1^- \oplus H_2^+ \oplus H_2^- $, we set
\[
  (\mathbf{q}, \mathbf{g}) = 
  \sum_{\alpha=1}^{2} ( (q^+_\alpha, g^+_\alpha)_\alpha^+ + 
  (q^-_\alpha, g^-_\alpha)_\alpha^-) ,
  \quad \|\mathbf{q}\| =  (\mathbf{q}, \mathbf{q})^{1/2} .
\] 

Now we construct the discrete analogs of differential operators
$\mathcal{D}_{\alpha}, \ \mathcal{D}^*_{\alpha}, \ \alpha = 1,2$
introduced according to (\ref{2.8}), (\ref{2.9}).
Using the above difference derivatives in space, we set
\begin{equation}\label{3.7}
  D_\alpha^+ y = - y_{x_\alpha},
  \quad \bm{x} \in \omega_\alpha^+ ,
  \quad  \alpha = 1,2 ,
\end{equation}
so that $D_\alpha^+: H \rightarrow H_\alpha^+, \ \alpha = 1,2$.
Similarly, we define $D_\alpha^-: H \rightarrow H_\alpha^-, \ \alpha = 1,2$, where
\begin{equation}\label{3.8}
  D_\alpha^- y = - y_{\bar{x}_\alpha},
  \quad \bm{x} \in \omega_\alpha^- ,
  \quad  \alpha = 1,2 .
\end{equation}
Thus
\begin{equation}\label{3.9}
  D: H \rightarrow V,
  \quad D = (D_1^+, D_1^-, D_2^+, D_2^- )^T .
\end{equation}
For the adjoint operator, we have
\begin{equation}\label{3.10}
  D^*: V \rightarrow H,
  \quad D^* = ((D_1^+)^*, (D_1^-)^*, (D_2^+)^*, (D_2^-)^* ),
\end{equation}
and
\begin{equation}\label{3.11}
  (D_\alpha^+)^* : H_\alpha^+ \rightarrow H ,
  \quad    (D_\alpha^+) q = q_{\bar{x}_\alpha},
\end{equation} 
\begin{equation}\label{3.12}
  (D_\alpha^-)^* : H_\alpha^- \rightarrow H ,
  \quad    (D_\alpha^-) q = q_{x_\alpha},
  \quad \bm{x} \in \omega,
  \quad  \alpha = 1,2 .
\end{equation}
The above discrete operators approximate the corresponding differential operators with the first order:
\begin{equation}\label{3.13}
  D_{\alpha}^\pm  u =  \mathcal{D}_{\alpha} u + \mathcal{O} (h_\alpha),
  \quad (D_{\alpha}^\pm)^* u =  \mathcal{D}_{\alpha}^* u + \mathcal{O} (h_\alpha),
  \quad  \alpha = 1,2 .
\end{equation}

For the operator-differential equation (\ref{2.13}), we put into the correspondence the equation
\begin{equation}\label{3.14}
  \frac{d y}{d t} + D^* K D y = \varphi (t),
  \quad 0 < y \leq T,  
\end{equation}
where, e.g, $\varphi (t) = f(\bm{x}, t), \  \bm{x} \in \omega$.
For equation (\ref{3.14}), we consider the Cauchy problem
\begin{equation}\label{3.15}
  y(0) = u^0 .
\end{equation} 

The construction of the operator $K$ is associated with the approximations (\ref{3.1})--(\ref{3.5}).
The most important properties are self-adjointness and positive
definiteness of the operator $K$.
The equation (\ref{3.14}) approximates the differential equation
(\ref{2.13}) with the second order.

The system of equations (\ref{2.10}), (\ref{2.11}) is attributed to the system
\begin{equation}\label{3.21}
  \frac{d y}{d t} + D^* \mathbf{g}  = \varphi (t),
  \quad 0 < t \leq  T ,
\end{equation} 
\begin{equation}\label{3.22}
  \mathbf{g}  = K D y .
\end{equation} 
For the flux problem (\ref{2.14}), (\ref{2.15}), we put into the correspondence the problem
\begin{equation}\label{3.23}
  C \frac{d \mathbf{g} }{d t} + D D^* \mathbf{g} = 
  D \varphi (t),
  \quad C = K^{-1},
  \quad 0 < t \leq  T ,
\end{equation} 
\begin{equation}\label{3.24}
  \mathbf{g}(0) = K D  u^0 .
\end{equation} 

Similarly to (\ref{2.20}), we prove the following estimate for the solution of the problem
(\ref{3.14}), (\ref{3.15}):
\begin{equation}\label{3.25}
  \|y(t)\| \leq \|u^0\| + \int_{0}^{t}  \|\varphi(\theta)\| d \theta .
\end{equation}
For the estimate (\ref{2.22}), we put into the correspondence the estimate
\begin{equation}\label{3.26}
  \|\mathbf{g}(t)\|_C \leq \|D u^0\|_K + \int_{0}^{t}  \|D \varphi(\theta)\|_K \, d \theta  
\end{equation} 
for the solution of the problem (\ref{3.23}), (\ref{3.24}).

\section{Operator-difference schemes}

We introduce a uniform grid in time with a step $\tau$  and let $y^n = y(t^n), \ t^n = n \tau$,
$n = 0,1, ..., N, \ N\tau = T$.
For numerical solving the problem (\ref{3.14}), (\ref{3.15}), we apply
the standard two-level scheme with weights, where equation (\ref{3.14})
is approximated by the scheme
\begin{equation}\label{4.1}
  \frac{y^{n+1} - y^{n}}{\tau }
  + A(\sigma y^{n+1} + (1-\sigma) y^{n}) = \varphi^n,
  \quad n = 0,1, ..., N-1,
\end{equation}
where
\begin{equation}\label{4.2}
  A = D^* K D,
  \quad A = A^* > 0  
\end{equation}
and, e.g., $\varphi^n = f(\sigma t^{n+1} + (1-\sigma) t^{n})$.
Taking into account (\ref{3.15}), the operator-difference equation (\ref{4.1})
is supplemented with the initial condition
\begin{equation}\label{4.3}
  y^0 = u^0 .
\end{equation}
The truncation error of the difference scheme (\ref{4.1})--(\ref{4.3}) is
$\mathcal{O} (|h|^2 + \tau^2 + (\sigma - 0.5) \tau)$. 

The study of the difference scheme is conducted
using the general theory of stability (well-posedness) for operator-difference
schemes \cite{Samarskii1989,SamarskiiMatusVabischevich2002}.
Let us formulate a typical result on stability
of difference schemes with weights for an evolutionary equation of first order.

{\bf Theorem  1.} {\it 
The scheme (\ref{4.1})--(\ref{4.3})
is unconditionally stable for  $\sigma \geq 0.5$,
and the difference solution satisfies the levelwise estimate
\begin{equation}\label{4.4}
  \|y^{n+1}\| \leq \|y^{n}\| + \tau \|\varphi^n\|,
  \quad n = 0,1, ..., N-1 .
\end{equation}
}

From (\ref{4.4}), in the standard way, we get the desired stability estimate
\[
  \|y^{n+1}\| \leq \|u^0\| + \sum_{k=0}^{n} \tau \|\varphi^k\| ,
\]
which may be treated as a direct discrete analogue of 
the a priori estimate (\ref{2.20}) for the solution of the differential problem
(\ref{2.12}), (\ref{2.18}).

Schemes with weights for a system of semi-discrete equations
(\ref{3.21}), (\ref{3.22}) are constructed in a similar way. We put
\begin{equation}\label{4.7}
  \frac{y^{n+1} - y^{n}}{\tau }
  + D^*(\sigma \mathbf{g}^{n+1} + (1-\sigma) \mathbf{g}^{n}) = \varphi^n,
  \quad n = 0,1, ..., N-1,
\end{equation}
\begin{equation}\label{4.8}
  \mathbf{g}^{n} = K D y^{n},
  \quad n = 0,1, ..., N .
\end{equation}
The scheme (\ref{4.7}), (\ref{4.7}) is equivalent to the scheme (\ref{4.1}).
In view of Theorem~1, it is stable under the restriction $\sigma \geq 0.5$,
and for the solution of difference problem (\ref{4.2}), (\ref{4.7}), (\ref{4.7}),
the a priori estimate (\ref{4.4}) holds.

The special consideration should be given to the flux problem (\ref{3.23}), (\ref{3.24}).
To solve it numerically, we apply the scheme
\begin{equation}\label{4.9}
  C \frac{\mathbf{g}^{n+1} - \mathbf{g}^{n}}{\tau }
  + D D^*(\sigma \mathbf{g}^{n+1} + (1-\sigma) \mathbf{g}^{n}) = D\varphi^n,
  \quad n = 0,1, ..., N-1,
\end{equation}
\begin{equation}\label{4.10}
  \mathbf{g}^0 = KD u^0 .
\end{equation} 

{\bf Theorem  2.} {\it 
The difference scheme (\ref{4.9}), (\ref{4.10})
is unconditionally stable for $\sigma \geq 0.5$,
and the difference solution satisfies the estimate
\begin{equation}\label{4.11}
  \|\mathbf{g}^{n+1}\|_C \leq \|\mathbf{g}^{n}\|_C + \tau \|D \varphi^n\|_K,
  \quad n = 0,1, ..., N-1 .
\end{equation}
}

From (\ref{4.11}), it follows the estimate
\[
  \|\mathbf{g}^{n+1}\|_C \leq \|D u^0\|_K + 
  \sum_{k=0}^{n} \tau \|D \varphi^k\|_K,
  \quad n = 0,1, ..., N-1 , 
\]
which corresponds to the estimate  (\ref{3.26}) for the solution of the problem
(\ref{3.23}), (\ref{3.24}).

The computational implementation of the unconditionally stable operator-difference schemes
(\ref{4.1})--(\ref{4.3}) for the parabolic equation
(\ref{2.1}) with mixed derivatives is
based on solving discrete elliptic problems at every time step.
For the problem (\ref{3.14}), (\ref{3.15}),
it seems more convenient to employ additive schemes (operator-splitting schemes)
that provide the transition to a new time level using
simpler problems associated with the inversion of the individual operators
$D^*_{\alpha} D_{\alpha},  \ \alpha =1,2$ rather then their combinations.
By the nature of the operators $D^*_{\alpha}, D_{\alpha},  \ \alpha =1,2$,
in this case, we speak of locally one-dimensional schemes.

The issues of designing unconditionally stable locally one-dimensional schemes 
for a parabolic equation without mixed derivatives have been studied in detain.
For parabolic equations with mixed derivatives, locally one-dimensional
schemes were constructed in several papers (see, e.g.,
\cite{mckee1970alternating,mishra2011stability}).
Strong results on unconditional stability
of operator-splitting schemes can be proved only in
a uninteresting case with
pairwise commutative operators (the equation with constant coefficients).
For our problems (\ref{2.1})--(\ref{2.4}),
the construction of locally one-dimensional schemes requires
separate consideration.

Let us investigate approaches to constructing locally one-dimensional schemes
for the problem (\ref{3.23}), (\ref{3.24}).
The computational implementation of the scheme with weights (\ref{4.9}), (\ref{4.10}),
which is unconditionally stable for $\sigma \geq 0.5$, is
associated with solving the system of difference equations for four components
of the vector $\mathbf{g}^{n+1}$.
The equations of this system are strongly coupled to each other, and this interconnection 
does exist not only for the spatial derivatives 
(operators $D^\pm_1 D^{\pm*}_2$, $D^\pm_2 D^{\pm*}_1$),
but also for the time derivatives ($k_{12} = k_{21} \neq 0$). 
Thus, we need to resolve the problem of splitting
for the operator at the time derivative, too.

The simplest case is splitting of the spatial operator without
coupling the time derivatives.
Such a technique is directly applicable for the construction
of locally one-dimensional schemes for parabolic equations
without mixed derivatives, where
\begin{equation}\label{5.1}
  k_{\alpha\beta} ({\bm x}) =k_{\beta\alpha}({\bm x})  = 0, 
  \quad \alpha \neq \beta = 1,2, ..., m,
  \quad {\bm x} \in \Omega  
\end{equation} 
in equation (\ref{2.1}).

Assume that
\[
 R = D D^*,
 \quad Q = \mathrm{diag} (D_1^+(D_1^+)^* , D_1^-(D_1^-)^* ,
 D_2^+(D_2^+)^* , D_2^-(D_2^-)^*) ,
\] 
i.e., $Q$ is the diagonal part of $R$.
For numerical solving the problem (\ref{3.23}), (\ref{3.24}),
we employ the difference scheme, where only the diagonal part of $R$
is shifted to the upper time level.
In our notation, we set
\begin{equation}\label{5.2}
  C \frac{\mathbf{g}^{n+1} - \mathbf{g}^{n}}{\tau }
  + Q(\sigma \mathbf{g}^{n+1} + (1-\sigma) \mathbf{g}^{n}) 
  + (R - Q)\mathbf{g}^{n} = D\varphi^n,
  \quad n = 0,1, ..., N-1,
\end{equation} 
with the initial conditions according to (\ref{4.10}). 

{\bf Theorem  3.} {\it 
The difference scheme (\ref{4.10}), (\ref{5.2})
is unconditionally stable for $\sigma \geq 2$,
and the difference solution satisfies the estimate 
\begin{equation}\label{5.3}
  \|\mathbf{g}^{n+1}\|_B \leq \|\mathbf{g}^{n}\|_B + \tau \|D \varphi^n\|_{B^{-1}},
  \quad n = 0,1, ..., N-1 ,
\end{equation}
where
\[
 B = C + \sigma\tau P - \frac{\tau}{2} R .
\] 
}

The scheme (\ref{4.10}), (\ref{5.2}) has the first-order approximation in time.
It seems more preferable, in terms of accuracy, to apply the scheme that is based on
the triangular decomposition of the self-adjoint matrix operator $\bm{R}$:
\begin{equation}\label{5.5}
  \bm{R} = \bm{R}_1 + \bm{R}_2,
  \quad  \bm{R}^*_1 = \bm{R}_2 .
\end{equation}
For the problem (\ref{3.23}), (\ref{3.24}), we construct the additive scheme
with the splitting (\ref{5.5}), where
\begin{equation}\label{5.7}
  (C + \sigma\tau R_1) C^{-1} (C + \sigma\tau R_2)  
  \frac{\mathbf{g}^{n+1} - \mathbf{g}^{n}}{\tau }
  + R \mathbf{g}^{n} = D\varphi^n,
  \quad n = 0,1, ..., N-1 .
\end{equation}
The main result is formulated in the following statement.

{\bf Theorem  4.} {\it 
The difference scheme (\ref{4.10}), (\ref{5.5})--(\ref{5.7})
is unconditionally stable for $\sigma \geq 0.5$,
and the difference solution satisfies the estimate (\ref{5.3}) with
\[
 B = (C + \sigma\tau R_1) C^{-1} (C + \sigma\tau R_2) - \frac{\tau}{2} R .
\] 
}

The alternating triangle operator-difference scheme
(\ref{4.10}), (\ref{5.5})--(\ref{5.7}) belongs to the class of schemes that are
based on a pseudo-time evolution process
---  the solution of the steady-state problem is obtained as a limit of this
pseudo-time evolution.
It has the second-order accuracy in time if $\sigma = 0.5$, and
ony the first order for other values of $\sigma $.

\end{document}